\documentclass[12pt]{article}

\usepackage{graphicx}
\usepackage{geometry}
\usepackage{amssymb}
\usepackage{amsmath,bm}
\usepackage{bbm}

\newcommand{\bey}{\begin{eqnarray}}
\newcommand{\eey}{\end{eqnarray}}
\newcommand{\beq}{\begin{equation}}
\newcommand{\eeq}{\end{equation}}

\newcommand{\sgn}{\text{sgn}}

\newtheorem{thm}{Theorem}[section]

\newtheorem{coro}{Corollary}[section]

\begin{document}

\vspace*{0in}

\noindent {\large \bf Estimable group effects for strongly correlated variables in linear models}

\bigskip
\noindent {MIN TSAO} \\ {\em Department of Mathematics \& Statistics, University of Victoria}

\bigskip

\bigskip

\bigskip


\bigskip

{

\noindent {\bf Abstract:} When a linear model has a multicollinearity problem caused by a group of strongly correlated predictor variables, it is well known that parameters for variables in such a group cannot be accurately estimated. We study linear combinations of these parameters to look for ones that can be accurately estimated. Under a uniform model for which the level of multicollinearity can be quantified, we find such linear combinations in a neighbourhood of a simple variability weighted average of these parameters. This variability weighted average actually benefits from multicollinearity in that the variance of its minimum-variance unbiased linear estimator is a monotone decreasing function of the level of multicollinearity, so it can be more accurately estimated at higher levels of multicollinearity. It is the only linear combination with this property. It can be substantially more accurately estimated than the average of parameters of uncorrelated variables. We also demonstrate that it can be accurately estimated outside the uniform model. Accurately estimated linear combinations have applications in inference and estimation for parameters of strongly correlated variables. As an alternative to the Ridge regression, we propose a method which uses the variability weighted average to estimate such parameters while retaining the least-squares estimates for parameters of other variables in the model.

\bigskip

}

\newpage

\section{Introduction}

Consider a linear model
\begin{equation}
\mathbf{y}=\mathbf{X}\bm{\beta} +\bm{\varepsilon}, \label{m0}
\end{equation}
where $\mathbf{y}$ is an $n\times 1$ vector of observations, $\mathbf{X}=[\mathbf{1},\mathbf{x}_1,\dots,\mathbf{x}_{q-1}]$ a
known $n\times q$ design matrix, $\boldsymbol{\beta}=(\beta_0,\beta_1,\dots,\beta_{q-1})^T$ an unknown $q \times 1$ vector of regression parameters, and $\boldsymbol{\varepsilon}$ an $n\times 1$ vector of random errors with mean zero and variance matrix $\sigma^2\mathbf{I}$.

Suppose the first
$p$ predictor variables $\mathbf{x}_1,\mathbf{x}_2,\dots,\mathbf{x}_p$ are strongly correlated ($2\leq p<q-1$). Then, this group of $p$ variables generates a multicollinearity problem for model (\ref{m0}). This type of multicollinearity problem arises often from data sets in observational studies where strongly correlated variables are not uncommon.
It manifests numerically through unusually large variances of least-squares estimators for parameters of the $p$ variables, sometimes accompanied by estimated parameters with unusually large absolute values. 
Reasons behind these numerical observations are well documented. They can be found in many books on linear models, {\em e.g.}, Draper and Smith (1998) and Belsley, Kuh and Welsch (2004). In particular, the latter book gives a comprehensive coverage on the detection of and remedies for the multicollinearity problem.

Although parameters of the strongly correlated variables $\beta_1$, $\beta_2$, $\dots,$ $\beta_p$ cannot be accurately estimated, other parameters in the model and some linear combinations of all $q$ parameters can still be. For convenience, we refer to a linear combination of $\beta_0, \beta_1, \dots, \beta_{q-1}$ as an {\em effect} of the underlying variables, and an effect is said to be {\em estimable} if it has an unbiased linear estimator whose variance is smaller than or comparable to the error variance $\sigma^2$, similar in magnitude to variances of least-squares estimators for parameters not affected by the multicollinearity problem.
Silvey (1969) gave the following result describing those effects that may be and those that may not be accurately estimated.
For simplicity, suppose $\mathbf{X}^T\mathbf{X}$ is nonsingular so that the least-squares estimator $\hat{\boldsymbol{\beta}}$ is available. Let $\mathbf{V}$ be the $q \times q$ orthogonal matrix whose columns $\mathbf{v}_i$ are orthonormal eigenvectors of $\mathbf{X}^T\mathbf{X}$, and let $\lambda_1\geq \lambda_2\geq \dots \geq \lambda_q>0$ be the eigenvalues of $\mathbf{X}^T\mathbf{X}$. For any fixed $q\times 1$ vector $\mathbf{c}$, there is a unique $q\times 1$ vector $\boldsymbol{\alpha}$ whose elements satisfy
\beq
\mathbf{c}=\alpha_1\mathbf{v}_1+\alpha_2\mathbf{v}_2+\dots+\alpha_q\mathbf{v}_q. \label{t1}
\eeq
Vector $\mathbf{c}$ defines an effect $\mathbf{c}^T\boldsymbol{\beta}$.
By the Guass-Markov theorem, the minimum-variance unbiased linear estimator for this effect is $\mathbf{c}^T\hat{\boldsymbol{\beta}}$. Using equation (\ref{t1}), the variance of this estimator $var(\mathbf{c}^T\hat{\boldsymbol{\beta}})$ can be shown to satisfy
\beq
var(\mathbf{c}^T\hat{\boldsymbol{\beta}})/\sigma^2=\frac{\alpha^2_1}{\lambda_1}+\frac{\alpha^2_2}{\lambda_2}+\dots+\frac{\alpha^2_q}{\lambda_q}.
\label{silvey}
\eeq

Under the constraint $\mathbf{c}^T\mathbf{c}=1$ which implies $\sum^p_{i=1}\alpha^2=1$, Silvey (1969) observed from equation (\ref{silvey}) that {``relatively precise estimation is possible in the directions of latent vectors of $\mathbf{X}^T\mathbf{X}$ corresponding to large latent roots; relatively imprecise estimation in these directions corresponding to small latent roots.''}
While one can use this observation to find $\mathbf{c}$ values such that $\mathbf{c}^T\hat{\boldsymbol{\beta}}$ have relatively small variances, it does not offer a meaningful interpretation for the underlying effects $\mathbf{c}^T{\boldsymbol{\beta}}$. Indeed, as Belsley, Kuh and Welsch (2004, p178) had noted that effects defined through the eigenvalues and eigenvectors of $\mathbf{X}^T\mathbf{X}$ are unlikely to be of practical interest. 

In this paper, we focus on a class of effects of practical interest and look among them for estimable effects. This class consists of all effects involving only $\beta_1, \beta_2, \dots, \beta_p$, so they have a clear interpretation as group effects for variables in the strongly correlated group. We are interested in finding estimable effects in this class; as none of the underlying parameters is estimable, such estimable effects are of theoretical interest and they also have practical applications some of which are mentioned in our motivations below. Specifically, we study the class of effects, $\Xi'$, given by
\beq
\Xi'=\{\xi(\mathbf{w'})\hspace{0.05in}|\hspace{0.05in}\xi(\mathbf{w'})=w'_1\beta_1+w'_2\beta_2+\dots+w'_p\beta_p \}, \label{lcombination}
\eeq
where $\mathbf{w'}=(w'_1,w'_2,\dots,w'_p)^T$ is any $p\times 1$ vector satisfying $\sum_{i=1}^p |w'_i|=1$. We call $\Xi'$ the class of {\em normalized group effects} of the $p$ strongly correlated variables. We choose constraint $\sum_{i=1}^p |w'_i|=1$ instead of $\sum_{i=1}^p (w'_i)^2=1$ because it allows $\Xi'$ to include commonly used weighted averages such as $\frac{1}{p}\sum^p_{i=1}\beta_i$, even though it is technically more difficult to handle as it is non-smooth. We call an effect $\xi(\mathbf{w'})$ in $\Xi'$ a {\em group effect} and the corresponding vector $\mathbf{w'}$ a {\em weight vector}. Individual parameters $\beta_1, \beta_2, \dots, \beta_p$ are special group effects in $\Xi'$ but they are not estimable. Our objective is to look for a characterization of the estimable effects in this class.

We are motivated by the fact that estimable effects are useful for inference and estimation concerning the underlying parameters $\beta_1, \beta_2, \dots, \beta_p$. They are also useful for knowing when accurate predictions can be made using the estimated model. For example, if the estimated value $\hat{\xi}(\mathbf{w'})$ of an estimable effect $\xi(\mathbf{w'})$ is significantly different zero, then we may reject the null hypothesis $H_0: \hspace{0.01in} \beta_1=\beta_2=\dots=\beta_p=0$ and conclude that one or more of the variables in this group are not zero. Also, each effect represents a linear constraint on the parameters, so an estimable effect can be used for dimension reduction. The parameter space for the underlying parameters is $\mathbb{R}^p$. If $\xi(\mathbf{w'})$ is estimable, then the unknown parameters satisfy
$$\hat{\xi}(\mathbf{w'})\approx w'_1\beta_1+w'_2\beta_2+\dots+w'_p\beta_p.$$
This reduces the parameter space from $\mathbb{R}^p$ to essentially a line in $\mathbb{R}^p$, and such a dimension reduction can be used for estimating the underlying parameters. The set of weight vectors $\{\mathbf{w}'\}$ associated with the estimable effects also defines a region in the space of the $p$ strongly correlated variables over which accurate predictions can be made using the estimated parameters $\hat{\beta}_1, \hat{\beta}_2, \dots, \hat{\beta}_p$, even though these are not accurate estimates of $\beta_1, \beta_2, \dots, \beta_p$.

To make the characterization problem manageable, we focus on a uniform model for which the level of multicollinearity generated by strongly correlated variables can be quantified and the impact of multicollinearity on group effects is mathematically tractable. Under this model, we find an optimal effect that benefits from multicollinearity in the sense that the variance of its minimum-variance unbiased linear estimator actually decreases when the level of multicollinearity increases; that is, this effect can be more accurately estimated when the level of multicollinearity goes higher. It is rather surprising that such a linear combination of the parameters exists as these parameters themselves are not estimable at high levels of multicollinearity. This optimal effect is the only one in class $\Xi'$ that benefits from multicollinearity; other effects all suffer in that the variances of their minimum-variance unbiased linear estimators all go to infinity as the level of multicollinearity approaches the extreme. The optimal effect has a simple interpretation as a variability weighted average of the underlying parameters. At any given level of multicollinearity,
all estimable effects are located around this effect. The uniform model captures effectively the impact of multicollinearity on group effects of strongly correlated variables in (\ref{m0}). We demonstrate through numerical examples that the variability weighted group effect for (\ref{m0}) is also estimable. As an example of its applications, we also discuss a constrained local regression method that uses the variability weighted average effect to estimate the underlying parameters. This method complements the Ridge regression (Horel and Kennard, 1970) and other penalized methods such as Lasso (see, {\em e.g.,} Hastie, Tibshirani and Wainwright, 2015) in that it is a local method for estimating the parameters of only strongly correlated variables;
the least-squares estimates of other parameters are unchanged.

The rest of this paper is organized as follows. In Section 2, we introduce the uniform model under which we reduce $\Xi'$ to a subclass $\Xi$. This subclass is only ``$(1/2^{p})$th'' the size of $\Xi'$ but it contains all effects that can be most accurately estimated. We then find the optimal variability weighted average effect through $\Xi$ and give a characterization of all estimable effects under the uniform model using this effect as a reference point. In Section 3, we present a numerical study on the variability weighted group effect and discuss the constrained local regression method. Proofs of theorems and corollaries are given in the Appendix.
It should be noted that strong correlations among predictor variables represents only one type of multicollinearity. We focus on this one type because it is the most common, and because it is local in nature so it can be isolated and modelled. The local nature of this type of multicollinearity is also discussed in the Appendix.

\section{The optimal and estimable effects in uniform models}
\subsection{The uniform model}

Consider a simple case of (\ref{m0})
\begin{equation}
\mathbf{y}=\mathbf{X}\boldsymbol{\beta} +\boldsymbol{\varepsilon} \label{m1}
\end{equation}
with $p$ predictors variables in $\mathbf{X}=[\mathbf{x}_1,\mathbf{x}_2,\dots,\mathbf{x}_p]$ where
each column $\mathbf{x}_i$ has a mean of $\bar{x}_i=0$ and a length of $\|\mathbf{x}_i\|=1$, and $\sigma^2=1$.
In the context of studying multicollinearity, there must be two or more variables in the model, so $p\geq 2$.
For this model, $\mathbf{X}^T\mathbf{X}$ is the correlation matrix; that is,
\beq
\mathbf{X}^T\mathbf{X} =
\left( \begin{array}{ccccc}
1 & r_{12} & r_{13} & \dots & r_{1p} \\
r_{21} & 1 & r_{23}  & \dots & r_{2p}  \\
r_{31} & r_{32} & 1 & \dots & r_{3p}  \\
\cdot&\cdot&\cdot&\dots &\cdot \\
\cdot&\cdot&\cdot&\dots &\cdot \\
r_{p1} & r_{p2}  & r_{p3} & \dots & 1 \\
\end{array} \right)_{p\times p},  \label{xtx}
\eeq
where $r_{ij}=corr(\mathbf{x}_i, \mathbf{x}_j)$. We assume that $n>p$ and variables $\mathbf{x}_i$ are linearly independent so that the least-squares estimator for
$\boldsymbol{\beta}$,
\beq
\hat{\boldsymbol{\beta}}= (\mathbf{X}^T\mathbf{X})^{-1}\mathbf{X}^T\mathbf{y}, \label{ole}
\eeq
is available. We define a {\em uniform model} as a linear model (\ref{m1}) that satisfies the following two
conditions:
\bey
& & \mbox{All $\mathbf{x}_i$ are positively correlated, {\em i.e.}, $r_{ij}>0$ for all $(i,j)$  .} \label{c1} \\
& & \mbox{All $\mathbf{x}_i$ are equally correlated, {\em i.e.}, $r_{ij}=r$ for some constant $r$.} \label{c2}
\eey
The uniform model is an approximation to a linear model with a group of strongly correlated predictor variables. If the group does not satisfy condition (\ref{c1}), we may change the signs of some of the variables so that it does without affecting the underlying model; the impact of this sign change is absorbed by a corresponding sign change of the parameters. Under condition
(\ref{c1}), as correlation coefficients of strongly correlated variables, $r_{ij}$ must all be close to but below 1. As such they are also close to each other, which implies condition (\ref{c2}) holds approximately. Theorem \ref{thm11} below gives a sufficient condition for obtaining (\ref{c1}) through sign changes of some variables in the group.

\begin{thm} \label{thm11}

 Let $\mathbf{x}_1, \mathbf{x}_2, \dots, \mathbf{x}_p $ be predictor variables in model (\ref{m1}) such that $\bar{x}_i=0$ and
 $\|\mathbf{x}_i\|=1$ for $i=1,2,\dots,p$. Suppose
 \beq
 |corr(\mathbf{x}_i, \mathbf{x}_1)|>\frac{\sqrt{2}}{2}, \hspace{0.3in} \mbox{for $i=2,3,\dots,p$}.\label{cond1}
\eeq
Then, among the $2^p$ sets of size $p$ each, formed by choosing exactly one element from
$\{\mathbf{x}_j, -\mathbf{x}_j\}$ for $j=1,2,\dots,p$, there exists a set from which every pair of variables has
a positive correlation.
\end{thm}

We call the set identified by Theorem \ref{thm11} an {\em all positive correlations arrangement} of the variables or an {\em APC arrangement} for short. Any group of variables satisfying condition (\ref{cond1}) has an APC arrangement that may be found by first computing $r_{1j}=corr(\mathbf{x}_1,\mathbf{x}_j)$ for $j=2,3, \dots, p$ and then changing all $\mathbf{x}_j$ where $corr(\mathbf{x}_1,\mathbf{x}_j)<0$ into $-\mathbf{x}_j$. The resulting APC arrangement is
 \[\mathbf{x}_1, \sgn(r_{12})\mathbf{x}_2, \dots, \sgn(r_{1p})\mathbf{x}_p.\]
Note that condition (\ref{cond1}) requires only the existence of one variable whose correlations with other variables are at least $\sqrt{2}/2$ (roughly 0.7) in absolute value. A strongly correlated group of variables behind a multicollinearity problem usually satisfies this condition. Thus, we may assume its variables are in an APC arrangement or equivalently condition (\ref{c1}) holds without loss of generality.

We use the uniform model for the simplicity it brings to the characterization problem. Specifically, with condition (\ref{c1}) we will be able to reduce the complexity of the characterization problem by a factor of $2^p$. With condition (\ref{c2}), we can measure the level of multicollinearity by using $r$; the larger $r$ is, the closer it is to 1, the more ill-conditioned $\mathbf{X}^T\mathbf{X}$ in (\ref{xtx}) and (\ref{ole}) becomes, and thus the higher the level of multicollinearity. This allows us to control the level of multicollinearity easily through a single parameter $r$, which makes it convenient to study the effect of multicollinearity quantitatively. The condition that $\bar{x}_i=0$ and $\|\mathbf{x}_i\|=1$ will be relaxed in Section 2.3.

\subsection{The average group effect and its optimality}

Let $\Xi'$ be the class of normalized group effects for parameters $\beta_1,\beta_2,\dots,\beta_p$ in (\ref{m1}).
Consider the subclass $\Xi=\{\xi(\mathbf{w})\}$ consisting of all properly weighted averages where
\beq
\Xi=\{\xi(\mathbf{w})\hspace{0.05in}|\hspace{0.05in}\xi(\mathbf{w})=w_1\beta_1+w_2\beta_2+\dots+w_p\beta_p \}, \label{eff}
\eeq
and $\mathbf{w}=(w_1,w_2,\dots,w_p)^T$, $\sum^p_{i=1} w_i=1$ and $w_i\geq 0$. This subclass contains only ``$(1/2^p)$'' of the effects in $\Xi'$ in that every effect $\xi(\mathbf{w})$ in $\Xi$ corresponds to $2^p$ effects $\xi(\mathbf{w}')$ in $\Xi'$ where each of the $2^p$ $\mathbf{w}'$ is obtained by changing the signs of some elements of $\mathbf{w}$ to negative. The constraint for this class $\sum^p_{i=1} w_i=1$ and $w_i\geq 0$ is a simple smooth constraint, which makes it easy to work with.

The minimum-variance unbiased linear estimator of $\xi(\mathbf{w})$ is
\beq
\hat{\xi}(\mathbf{w})=\mathbf{w}^T\hat{\boldsymbol{\beta}}=w_1\hat{\beta}_1+w_2\hat{\beta}_2+\dots+w_p\hat{\beta}_p,
\label{eeff}
\eeq
where $\hat{\boldsymbol{\beta}}=(\hat{\beta}_1,\hat{\beta}_2,\dots,\hat{\beta}_p)^T$ is the least-squares
estimator of $\boldsymbol{\beta}$ in (\ref{ole}).
The optimal effect in $\Xi$ is defined as the one whose estimator $\hat{\xi}(\mathbf{w})$ has the smallest variance.
The average of the parameters is an effect in $\Xi$ defined by weight vector $\mathbf{w}_0=(1/p, 1/p, \dots, 1/p)^T\in \mathbb{R}^p.$ We denote this effect by $\tau_a$ and call it the {\em average group effect} of the variables; that is,
\beq
\tau_a=\xi(\mathbf{w}_0)=\mathbf{w}_0^T\boldsymbol{\beta}=\frac{1}{p} \sum^p_{i=1} \beta_i . \label{ge}
\eeq
Correspondingly, let $\mathbf{w}_j=(w_1,w_2,\dots,w_p)^T$ where $w_j=1$ and $w_k=0$ for $k\neq j$. We call $$\beta_j=\xi(\mathbf{w}_j)$$ the {\em individual effect} of variable $\mathbf{x}_j$. We will first show that the average group effect $\tau_a$ is the optimal in $\Xi$ and then show that it is also the optimal over the much larger class $\Xi'$. To this end, we present Theorem \ref{thm30} below which gives the variance of $\hat{\xi}(\mathbf{w})$.

\begin{thm} \label{thm30}

Suppose conditions (\ref{c1}) and (\ref{c2}) hold. Then, the variance of a minimum-variance unbiased linear estimator $\hat{\xi}(\mathbf{w})$ in (\ref{eeff}) is
\beq
var(\hat{\xi}(\mathbf{w}))=
\frac{[1+(p-2)r]\sum^p_{i=1}w_i^2 - 2r\sum_{i<j}w_iw_j   }
{1+(p-2)r-(p-1)r^2}.  \label{effvar0}
\eeq

\end{thm}

The following corollaries highlight the impact of multicollinearity on the variances of estimators for individual and average group effects. To emphasize their dependence on the level of multicollinearity $r$, we will write the variances as $var(\hat{\xi}(\mathbf{w}),r)$.
For the estimator of an individual effect $\hat{\beta}_j$, we have

\begin{coro} \label{coro2}
Suppose conditions (\ref{c1}) and (\ref{c2}) hold.
Then, for a fixed $r$,
\beq
var(\hat{\beta_j},r)=var(\hat{\xi}(\mathbf{w}_j),r)=\max_\mathbf{w}\{var(\hat{\xi}(\mathbf{w}),r)\}.  \label{re2}
\eeq
Further, $var(\hat{\beta}_j,r)$ is a strictly monotone increasing function of $r$ and
\beq
\lim_{r\rightarrow 1} var(\hat{\beta}_j,r) = + \infty.  \label{limit1}
\eeq
\end{coro}

Corollary \ref{coro2} confirms the known result that the variance of the unbiased estimator for $\beta_j$ goes to infinity when the level of multicollinearity goes to the extreme. It also shows that among the subclass of effects $\Xi$ defined in (\ref{eff}), individual parameters are the most difficult to estimate in that the variances of their unbiased estimators $\hat{\beta}_j=\hat{\xi}(\mathbf{w}_j)$ are the largest among unbiased estimators (\ref{eeff}) for this subclass of effects. This observation is also valid over the larger class $\Xi'$.
We now give the fixed-$r$ and asymptotic optimality of the average group effect $\tau_a$.
\begin{coro} \label{coro3}
Suppose conditions (\ref{c1}) and (\ref{c2}) hold.
Then, at any fixed $r$
\beq
var(\hat{\tau}_a,r)=var(\hat{\xi}(\mathbf{w}_0),r)=\min_\mathbf{w}\{var(\hat{\xi}(\mathbf{w}),r)\}.  \label{re1}
\eeq
Also, $var(\hat{\tau}_a,r)$ is a strictly monotone decreasing function of $r$ with
\beq
\lim_{r\rightarrow 1} var(\hat{\tau}_a,r) = \frac{1}{p^2}. \label{limit2}
\eeq
Further, among variances $var(\hat{\xi}(\mathbf{w}),r)$ of unbiased estimators (\ref{eeff}) for effects in $\Xi$, $var(\hat{\tau}_a,r)$ is the only one that remains bounded as $r\rightarrow 1$.
\end{coro}

Corollary \ref{coro3} shows the average group effect $\tau_a$ is the optimal effect in $\Xi$ as its minimum-variance unbiased estimator $\hat{\xi}(\mathbf{w}_0)$ has the smallest variance among all estimators (\ref{eeff}) for effects in $\Xi$. Its proof also shows at any fixed $r$, other estimable effects ${\xi}(\mathbf{w})$ are all around ${\xi}(\mathbf{w}_0)$ in that their associated weight vectors $\mathbf{w}$ are in a neighborhood of $\mathbf{w}_0$; see the remark after the proof of Corollary \ref{coro3} in the Appendix for a detailed discussion. The proof also shows that for any fixed $\mathbf{w}\neq \mathbf{w}_0$, $var(\hat{\xi}(\mathbf{w}),r)\rightarrow \infty$ as $r\rightarrow 1$.
Because the minimum value $var(\hat{\xi}(\mathbf{w}_0),r)$ is a decreasing function of $r$, the convergence to infinity by $var(\hat{\xi}(\mathbf{w}),r)$ for $\mathbf{w}\neq \mathbf{w}_0$ is not uniform. There is in fact a small neighborhood of $\mathbf{w}_0$, say $\{\mathbf{w} \hspace{0.05in} | \hspace{0.05in} \mathbf{w} \in \mathbb{R}^p, \hspace{0.05in}  \|\mathbf{w}-\mathbf{w}_0\| < \delta(r)\}$, over which $var(\hat{\xi}(\mathbf{w}),r)$ becomes smaller when $r$ increases. But this neighborhood becomes smaller and smaller as its radius $\delta(r) \rightarrow 0$ when $r\rightarrow 1$.

Table 1 contains $var(\hat{\tau}_a,r)$ and $var(\hat{\beta_j},r)$ values at various $r$ levels for
a uniform model with $p=8$ variables. There are two points to note. (i) The $r$ value in column 1 starts at $r=0$ which corresponds to an orthogonal design, and increases to $r=0.999$ which represents very strong correlation. In column 2, $var(\hat{\tau}_a,r)$ decreases from $0.125$ to $0.015$; the latter is about 12\% of the former. Thus the average group effect is much more accurately estimated under strong correlations among the variables. (ii) Comparing columns 2 and 3, as $r$ increases we see a clear divergence of the two variances with $var(\hat{\tau}_a,r)$ going down to its limit of
$1/p^2=1/8^2$ and $var(\hat{\beta_j},r)$ going up to infinity. This shows the type of multicollinearity generated by strong correlations among variables represents a redistribution of information in favor of the average group effect $\tau_a$ at the expense of individual effects $\beta_j$.

\begin{table}
\caption{\label{table1} Variances of the average group effect estimator $var(\hat{\tau}_a,r)$ and
individual effect estimator $var(\hat{\beta_j},r)$ at various levels of multicollinearity $r$.}
\centering
\fbox{%
\begin{tabular}{*{5}{c}}
\hline
$r$ & \hspace{0.3in} & $var(\hat{\tau}_a,r)$ & \hspace{0.3in} & $var(\hat{\beta_j},r)$ \\ \hline
0.0000000 & & 0.12500000 & & 1.000000 \\
0.5000000 & & 0.02777778 & & 1.777778\\
0.6666667 & & 0.02205882 & & 2.647059\\
0.7500000 & & 0.02000000 & & 3.520000\\
0.8000000 & & 0.01893939 & & 4.393939\\
0.8333333 & & 0.01829268 & & 5.268293\\
0.8571429 & & 0.01785714 & & 6.142857\\
0.8750000 & & 0.01754386 & & 7.017544\\
0.8888889 & & 0.01730769 & & 7.892308\\
0.9000000 & & 0.01712329 & & 8.767123\\
0.9990000 & & 0.01563868 & & 875.015639\\ \hline
\end{tabular}}
\end{table}

We now go back to the class of all normalized group effects $\Xi'$ in (\ref{lcombination}). For an effect $\xi(\mathbf{w'})$ in $\Xi'$ but not in $\Xi$, one or more of the elements of $\mathbf{w}'=(w'_1,w'_2,\dots,w'_p)^T$ must be negative. Consider the corresponding effect $\xi(\mathbf{w})$ in $\Xi$ where $\mathbf{w}=(w_1,w_2,\dots,w_p)^T$ and $w_i=|w_i'|$ for $i=1,2,\dots,p$. Noting that Theorem \ref{thm30} applies to minimum-variance unbiased estimators (\ref{eeff}) for all effects $\xi(\mathbf{w'})$ in $\Xi'$ under the uniform model and since
\[
\sum^p_{i=1}w_i^2=\sum^p_{i=1}(w'_i)^2 \mbox{\hspace{0.1in} but \hspace{0.1in}} \sum^p_{i<j}w_iw_j \geq \sum^p_{i<j}w'_iw'_j,
\]
by (\ref{effvar0}), we have
\beq
var(\hat{\xi}(\mathbf{w})) \leq var(\hat{\xi}(\mathbf{w'})). \label{compare}
\eeq
In the sense of (\ref{compare}), the subclass $\Xi$ contains all effects of $\Xi'$ that can be most accurately estimated, so the average group effect $\tau_a$ is also the optimal over the much larger class of effects $\Xi'$. For any $\mathbf{w}'\neq \mathbf{w}_0$, we can also prove that $var(\hat{\xi}(\mathbf{w'}),r)$ goes to infinity as $r$ approaches 1.

\subsection{The variability weighted group effect for a general uniform model}

Consider a more general linear model
\begin{equation}
\mathbf{y}=\beta_0\mathbf{1}_n+\mathbf{X}\boldsymbol{\beta} +\boldsymbol{\varepsilon}, \label{m2}
\end{equation}
where $\mathbf{y}$ is an $n\times 1$ vector of response variable values, $\beta_0$ the unknown intercept
term, $\mathbf{1}_n$ the $n\times 1$ vector of $1$'s, and
$\mathbf{X}=[\mathbf{x}_1,\mathbf{x}_2,\dots,\mathbf{x}_p]$ the known $n\times p$
design matrix,  $\boldsymbol{\beta}=(\beta_1,\beta_2,\dots,\beta_p)^T$ the unknown vector of
parameters of variables $\mathbf{x}_i$, and $\boldsymbol{\varepsilon}$ the $n\times 1$ vector of random error satisfying
\beq
E(\boldsymbol{\varepsilon})=0 \mbox{\hspace{0.1in} and \hspace{0.1in}}
var(\boldsymbol{\varepsilon})=\sigma^2\mathbf{I}. \label{eterm2}
\eeq
We assume variables in (\ref{m2}) satisfy $corr(\mathbf{x}_i,\mathbf{x}_j)=r>0$ for all $(i,j)$ and call such a model the {\em general uniform model}. Comparing with the simple uniform model (\ref{m1}), the general uniform model has an intercept and there are no restrictions on the mean or
length of variables $\mathbf{x}_i$. Let
\[ s_j=\left[\sum^n_{i=1} (x_{ij}-\bar{x}_j)^2\right]^{1/2}, \]
where $x_{ij}$ is the $i$th element of variable $\mathbf{x}_j$ and $\bar{x}_j$ the mean of its elements for
$j=1,2,\dots,p$. Without loss of generality, we assume that $s_j>0$ for all $j$.
To focus on the impact of the level of multicollinearity $r$ on effects in $\Xi'$, we hold $s_j$ which also affect the accuracy of estimation fixed, so our results here are conditional on $s_j$. Let
$\mathbf{w}_w=(w_1^*,w_2^*,\dots,w_p^*)^T$ where
\beq
w^*_j=\frac{s_j}{\sum^p_{i=1}s_i}. \label{weights}
\eeq
We define a {\em variability weighted group effect}, $\tau_w$, for variables in (\ref{m2}) as
\beq
\tau_w=\mathbf{w}^T_w\boldsymbol{\beta}=\sum^p_{i=1} w^*_i\beta_i.  \label{wef}
\eeq
In (\ref{wef}), the weight $w^*_i$ for parameter $\beta_i$ is large if the variability of $\mathbf{x}_i$,
represented by $s_i$, is large relative to that of other variables. Since parameters for variables with
larger variabilities tend to be more accurately
estimated, this variability weighted group effect $\tau_w$ has the appeal of giving more weight to the more
accurately estimated individual effects. More importantly, $\tau_w$ corresponds to the average group effect for
variables in the simple uniform model (\ref{m1}). This connection gives it the optimality enjoyed by
the average group effect in (\ref{m1}). Theorem \ref{thm3} below further
illustrates this point.

\begin{thm} \label{thm3}
For model (\ref{m2}), suppose $corr(\mathbf{x}_i,\mathbf{x}_j)=r>0$ for all $(i,j)$. Let
$\hat{\boldsymbol{\beta}}$ be the least-squares estimator of $\boldsymbol{\beta}$ and
$\hat{\tau}_w$ be the minimum-variance unbiased linear estimator for $\tau_w$ where
\beq
\hat{\tau}_w=\mathbf{w}^T_w\hat{\boldsymbol{\beta}}=\sum^p_{i=1} w^*_i\hat{\beta}_i.  \label{lsetauw}
\eeq
Then, $var(\hat{\tau}_w,r)$ is a strictly monotone decreasing function $r$ and
\beq
\lim_{r\rightarrow 1} var(\hat{\tau}_w,r) = \frac{\sigma^2}{(\sum^p_{i=1}s_i)^2}.  \label{limit5}
\eeq
Further, among all estimators (\ref{eeff}) for group effects in $\Xi'$ for model (\ref{m2}), $\hat{\tau}_w$ is the only one with a bounded variance
when $r$ approaches 1.
\end{thm}

By Theorem \ref{thm3}, $\tau_w$ is the optimal effect among the class of normalized group effects $\Xi'$ for model (\ref{m2}) in an asymptotic sense as the variance of its estimator $var(\hat{\xi}(\mathbf{w}_w),r)$ is asymptotically the smallest as $r$ approaches 1. For a fixed $r$, the variance of an estimator $var(\hat{\xi}(\mathbf{w}),r)$ also depends on $s_i$, so $var(\hat{\xi}(\mathbf{w}_w),r)$ may not be the smallest. Thus $\tau_w$ may not be the exact optimal but it is nearly optimal when $r$ is large. All estimable effects for model (\ref{m2}) are in a neighborhood of $\tau_w$.
The average group effect $\tau_a$ for (\ref{m2}) is estimable when $s_i$ do not differ too much. In this case, individual elements of $\mathbf{w}_w$ are all roughly $1/p$ and the average group effect is in the neighborhood of $\tau_w$. At a fixed $r$, the most difficult effect to estimate is the individual effect $\beta_j$ of the variable $\mathbf{x}_j$ with the smallest variablility $s_j$. See the proof of Theorem \ref{thm3} in the Appendix and the remark after the proof for details.

\section{Numerical examples and applications}

We have conducted a detailed analysis of the multicollinearity problem in uniform models (\ref{m1}) and (\ref{m2}), which are made of a single group of strongly correlated variables with a uniform correlation structure. For a general model (\ref{m0}) having one or more groups of strongly correlated variables with non-uniform correlation structures as well as not strongly correlated variables, a similarly detailed analysis of the problem seems to be difficult. Fortunately, the multicollinearity problem due to a group of strongly correlated variables is local in nature in that parameters for variables outside of the group are little affected by the problem. Also, as we have noted in the discussion following (\ref{c1}) and (\ref{c2}) that the uniform structure is a good approximation to the correlation structure of a group of strongly correlated variables in its APC arrangement. So the uniform models provide a useful local approximation to the general model (\ref{m0}) in terms of this type of multicollinearity problem.

Indeed, much of what we have learned through the uniform models are still useful for a general model (\ref{m0}). Specifically, (i) the within group APC arrangement is effective for finding estimable effects of a group of strongly correlated variables in model (\ref{m0}), (ii) the variability weighted group effect $\tau_w$ of such a group defined under its APC arrangement is estimable, in fact approximately optimal, and it can be substantially more accurately estimated than the average of parameters of the same number of uncorrelated variables, and (iii) $\tau_w$ still serves as a reliable location around which other estimable effects may be found. We now illustrate these observations with numerical examples. We also apply $\tau_w$ for inference and estimation of the underlying parameters.

\subsection{Numerical examples}

For economy of space, we use the following setup for all examples in this subsection. Consider a linear model (\ref{m0}) with 10 predictor variables $\mathbf{x}_1,\mathbf{x}_2,\dots,\mathbf{x}_{10}$ in three groups
$\{\mathbf{x}_1,\mathbf{x}_2\}$, $\{\mathbf{x}_{3}, \mathbf{x}_4,\mathbf{x}_5\}$ and $\{\mathbf{x}_6,\mathbf{x}_7,\dots,\mathbf{x}_{10}\}$. We use two parameters $w_1, w_2\in [0,1]$ and 10 independent $n$-variate standard normal random vectors $\mathbf{z}_i$ to generate different levels of within group correlations as follows:
\bey
& & \mathbf{x}_1=\mathbf{z}_1, \hspace{0.1in}  \mathbf{x}_2=w_1\mathbf{z}_1+(1-w_1)\mathbf{z}_2; \nonumber \\
& & \mathbf{x}_3=\mathbf{z}_3, \hspace{0.1in}  \mathbf{x}_4=w_1\mathbf{z}_3+(1-w_1)\mathbf{z}_4, \hspace{0.1in}
\mathbf{x}_5=w_2\mathbf{z}_3+(1-w_2)\mathbf{z}_5, \label{expl} \\
& & \mathbf{x}_j=\mathbf{z}_j, \mbox{ \hspace{0.00in} for $j=6,7,8,9,10$}. \nonumber
\eey
So this model (\ref{m0}) contains two separate groups of correlated predictor variables as well as five independent predictor variables. The theoretical non-zero correlation coefficients are:
\bey
& & \sigma_{12}=\sigma_{21}=w_1[w_1^2+(1-w_1)^2]^{-1/2} ; \nonumber \\
& &  \sigma_{34}=\sigma_{43}=w_1[w_1^2+(1-w_1)^2]^{-1/2}, \hspace{0.1in}
\sigma_{35}=\sigma_{53}=w_2[w_2^2+(1-w_2)^2]^{-1/2}, \nonumber\\
& & \sigma_{45}=\sigma_{54}=w_1w_2\{[w_1^2+(1-w_1)^2][w_2^2+(1-w_2)^2]\}^{-1/2}. \nonumber
\eey
These correlation coefficients are all positive, so the two groups of correlated variables come with the APC arrangement.
The numerical examples given below involve observed values of variables $\mathbf{x}_i$. There are mild correlations among the observed values of the last five variables. The observed correlation coefficients among variables in the two correlated groups also differ somewhat from the theoretical values above, and they do not have a uniform structure.

For all cases of this model (\ref{m0}) that we discuss below, the sample size and the true parameter values are:  $n=15$, $\sigma^2=1$ and
\beq \boldsymbol{\beta}=(\beta_0,\beta_1,\beta_2,\dots,\beta_{10})^T=(5,0,0,1,2,3,1,1,1,2,3)^T. \label{beta} \eeq
For ease of comparison, the 10 variables are conveniently organized in two pairs of groups, $\{\mathbf{x}_1,\mathbf{x}_2\}$ and $\{\mathbf{x}_{6}, \mathbf{x}_7\}$,  and $\{\mathbf{x}_{3}, \mathbf{x}_4,\mathbf{x}_5\}$ and $\{\mathbf{x}_8,\mathbf{x}_9,\mathbf{x}_{10}\}$. In each pair, both groups have the same number of variables but only the first group contains correlated variables. We compare estimated values and their variances of the following average and variability weighted group effects:

\vspace{0.1in}

\hspace*{0.0in} 1. $\tau_1$: the average group effect for correlated group $\{\mathbf{x}_1,\mathbf{x}_2\}$.\\
\hspace*{0.2in} 2. $\tau_2$: the average group effect for correlated group $\{\mathbf{x}_{3}, \mathbf{x}_4,\mathbf{x}_5\}$.\\
\hspace*{0.2in} 3. $\tau_3$: the average group effect for independent group $\{\mathbf{x}_{6}, \mathbf{x}_7\}$.\\
\hspace*{0.2in} 4. $\tau_4$: the average group effect for independent group $\{\mathbf{x}_8,\mathbf{x}_9,\mathbf{x}_{10}\}$.\\
\hspace*{0.2in} 5. $\tau_1^w$: the variability weighted group effect for correlated group $\{\mathbf{x}_1,\mathbf{x}_2\}$.\\
\hspace*{0.2in} 6. $\tau_2^w$: the variability weighted group effect for correlated group $\{\mathbf{x}_{3}, \mathbf{x}_4,\mathbf{x}_5\}$.\\
\hspace*{0.2in} 7. $\tau_3^w$: the variability weighted group effect for independent group $\{\mathbf{x}_{6}, \mathbf{x}_7\}$.\\
\hspace*{0.2in} 8. $\tau_4^w$: the variability weighted group effect for independent group $\{\mathbf{x}_8,\mathbf{x}_9,\mathbf{x}_{10}\}$.

\vspace{0.1in}

Table \ref{table2} to Table \ref{table6} contain numerical results for five different cases defined by different $(w_1,w_2)$ values. Two of these involve adjustments to the predictor variables $\mathbf{x}_i$ for illustrating the importance of the APC arrangement and the advantage of $\tau^w_i$ over $\tau_i$. Each table represents one case, and it contains the mean and variance of 1000 estimates (\ref{eeff}) of each $\tau$ in the list above. These 1000 estimates are generated as follows. We first generate a design matrix $\mathbf{X}$ through randomly chosen $\mathbf{z}_i$ values and (\ref{expl}) for the case. Then, with this same design matrix we generate 1000 $\mathbf{y}$-observations, $\mathbf{y}_1, \mathbf{y}_2,\dots,\mathbf{y}_{1000}$, each of length $n=15$, using parameters in (\ref{beta}). The resulting 1000 pairs of $(\mathbf{y}_i,\mathbf{X})$ are then used to compute 1000 least-squares estimates $\hat{\boldsymbol{\beta}}$ and then 1000 estimates (\ref{eeff}) for each $\tau$. For example, with a pair of $(\mathbf{y}_i,\mathbf{X})$ the average group effect $\tau_2$ and the weighted group effect $\tau^w_2$ are estimated by
\beq
 \hat{\tau}_2=\frac{1}{3}(\hat{\beta}_3+\hat{\beta}_4+\hat{\beta}_5) \mbox{\hspace{0.1in} and \hspace{0.1in}}
\hat{\tau}_2^w=w^*_1\hat{\beta}_3+w^*_2\hat{\beta}_4+w^*_3\hat{\beta}_5   \label{t5.1}
\eeq
where $\hat{\beta}_j$ are the least-squares estimates of $\beta_j$, and weights $w^*_j$ are defined by (\ref{weights}) and computed using the observed values of the second group of correlated variables $\{\mathbf{x}_{3}, \mathbf{x}_4,\mathbf{x}_5\}$.
As a reference, we also include in each table the means and variances of 1000 least-squares estimates for 8 individual effects $\beta_j$ (only 8 are included to fit the space). Although the numerical results are tied to the $\mathbf{X}$ used, since this simulation can be easily reproduced and results do not depend on $\mathbf{X}$  very much as long as it is generated randomly through (\ref{expl}), we do not show the $\mathbf{X}$ matrix used.

Case 1 (Table 2) represents a model containing two separate weakly correlated groups. With $(w_1,w_2)=(0.3,0.4)$, the observed correlations from the design matrix $\mathbf{X}$ used to compute Table 2 range from 0 to 0.57. So there is no strong correlation or multicollinearity here. The right-side of Table 2 shows the least-squares estimators for individual effects have small biases and variances. The left-side shows all 8 group effects are accurately estimated. Here, because the estimators such as (\ref{t5.1}) are unbiased, they are accurate when their variances are small. Comparing the variances of estimators for $\tau_i$ with that for $\tau^w_i$, we see the latter are slightly smaller in this case of weakly correlated data.

\begin{table}
\caption{\label{table2} Case 1: Mean and variance of 8 estimated group effects and 8 estimated individual effects based on 1000 simulated values. Model setting: $(w_1,w_2)=(0.3,0.4)$.}
\centering
\fbox{%
\begin{tabular}{*{6}{c}}
\hline
Effect &       Mean         & Variance  &  Effect  &    Mean     & Variance \\ \hline
$\tau_1$   & -0.02449476 & 0.20621566&$\beta_1$ &-0.03199310 & 0.3751411 \\
$\tau_2$   & 2.00295008  & 0.10476791&$\beta_2$ &-0.01699641 &0.7215614 \\
$\tau_3$   & 1.01228831  & 0.18901021&$\beta_3$ & 1.01033591 &0.3592244 \\
$\tau_4$   & 1.99792105  & 0.07997519&$\beta_4$ & 1.96807691 &1.0212519\\
$\tau_1^w$ & -0.02534364 & 0.19855636&$\beta_5$ & 3.03043743 &0.2736429 \\
$\tau_2^w$ & 1.97757125  & 0.08983408&$\beta_6$ & 1.01754264 &0.4063596 \\
$\tau_3^w$ & 1.01139891  & 0.17941193&$\beta_9$ & 2.00708677 &0.1447360\\
$\tau_4^w$ & 1.90246465  & 0.06614164&$\beta_{10}$ & 2.98164306 & 0.6757436 \\ \hline
\end{tabular}}
\end{table}

Case 2 (Table 3) represents a model containing two strongly correlated groups. With $(w_1,w_2)=(0.9,0.95)$, the observed correlations among variables within each of the first two groups are all over 0.90. So there is multicollinearity from these two groups/sources. The right-side of Table 3 shows the least-squares estimators for individual effects of the 5 variables in these two groups doing badly with large variances. But that for the 5 (theoretically) uncorrelated individual effects (only 3 are shown) are not affected and their variances are comparable to the corresponding variances in Table 2. The left-side shows all 8 group effects are also accurately estimated with small variances. In particular, variances for weighted group effects of the two correlated groups ($\tau^w_1$ and $\tau_2^w$) dropped substantially from Case 1, indicating that increased correlation/multicollinearity is making them more accurately estimated, whereas variances of other effects changed relatively little.

\begin{table}
\caption{\label{table3} Case 2: Mean and variance of 8 estimated group effects and 8 estimated individual effects based on 1000 simulated values. Model setting: $(w_1,w_2)=(0.90,0.95)$.}
\centering
\fbox{%
\begin{tabular}{*{6}{c}}
\hline
Effect &       Mean         & Variance  &  Effect  &    Mean     & Variance \\ \hline
$\tau_1$   & -0.024747684 & 0.11120434 &$\beta_1$ & 0.3589857  & 43.4826062 \\
$\tau_2$   & 2.008298266  & 0.14466301 &$\beta_2$ & -0.4084811 & 51.7201745 \\
$\tau_3$   & 0.980029467  & 0.16681453 &$\beta_3$ & 0.8285308  & 63.2663587 \\
$\tau_4$   & 1.993633063  & 0.05175908 &$\beta_4$ & 1.9995771  & 14.6978134 \\
$\tau_1^w$ & -0.008843552 & 0.03207452 &$\beta_5$ & 3.1967869  & 82.7437839 \\
$\tau_2^w$ & 1.976905206  & 0.05521688 &$\beta_6$ & 0.9990162  & 0.1689072 \\
$\tau_3^w$ & 0.982124495  & 0.15422274 &$\beta_9$ & 2.0013250  & 0.1080171 \\
$\tau_4^w$ & 1.959080262  & 0.05942050 &$\beta_{10}$ & 2.9805580 & 0.3040263 \\ \hline
\end{tabular}}
\end{table}

\begin{table}
\caption{\label{table4} Case 3: Mean and variance of 8 estimated group effects and 8 estimated individual effects based on 1000 simulated values. Model setting: $(w_1,w_2)=(0.999,0.999)$.}
\centering
\fbox{%
\begin{tabular}{*{6}{c}}
\hline
Effect &       Mean         & Variance  &  Effect  &    Mean     & Variance \\ \hline
$\tau_1$   & 0.009940493 & 0.06555451 &$\beta_1$ &    -29.1023643 & 2.226207e+05 \\
$\tau_2$   & 2.000834123 & 0.06990159 &$\beta_2$ &    29.1222453  & 2.229771e+05 \\
$\tau_3$   & 0.986230082 & 0.20504053 &$\beta_3$ &    1.8645331   & 4.262174e+05 \\
$\tau_4$   & 2.011184589 & 0.16574333 &$\beta_4$ &    11.4377567  & 3.997136e+05 \\
$\tau_1^w$ &-0.003200664 & 0.03051943 &$\beta_5$ &    -7.2997874  & 1.156558e+05 \\
$\tau_2^w$ & 2.000684397 & 0.01678193 &$\beta_6$ &     0.9772109  & 1.769559e-01\\
$\tau_3^w$ & 0.987087119 & 0.21823662 &$\beta_9$ &    2.0260632   & 4.473404e-01 \\
$\tau_4^w$ & 2.146899666 & 0.14498650 &$\beta_{10}$ & 2.9872200   & 3.417891e-01 \\ \hline
\end{tabular}}
\end{table}

Case 3 (Table 4) represents a model with extreme multicollinearity. With $(w_1,w_2)=(0.999,0.999)$, the observed correlations among variables within each of the first two groups are all over 0.999. The right-side of Table 5 shows the least-squares estimators for individual effects of the 5 correlated variables brake down with variances in the order of $10^5$. But that for the 5 uncorrelated effects (only 3 are shown) are again not affected, illustrating the local nature of this type of multicolinearity problem. The left-side shows all 8 group effects are also accurately estimated with small variances. In particular, variances for $\tau^w_1$ and $\tau^w_2$ are now so small that they are close to the theoretical lower limits (\ref{limit5}) for the uniform model (not shown). The small variances associated with the average group effects $\tau_1$ and $\tau_2$ of the two correlated groups reflect the remark after Theorem \ref{thm3}; the $s_i$ of the variables here are all close, so the average group effect is an estimable effect.
This observation is further confirmed in Case 4 (Table 5) where the level of multicollinearity has been lowered to $(w_1,w_2)=(0.99,0.99)$ from Case 3, but variables $\mathbf{x}_2$ and $\mathbf{x}_5$ are replaced with $2\mathbf{x}_2$ and $2\mathbf{x}_5$. This change creates a situation where $s_2$ dominates $s_1$, and $s_5$ dominates $s_3$ and $s_4$. So the average group effects $\tau_1$ and $\tau_2$ are no longer near the weighted group effects $\tau^w_1$ and $\tau^w_2$. Consequently, they are not estimable and their variances become very large whereas that of the two weighted group effects are not affected by this change.

\begin{table}
\caption{\label{table5} Case 4: Mean and variance of 8 estimated group effects and 8 estimated individual effects based on 1000 simulated values. Model setting: $(w_1,w_2)=(0.99,0.99)$ and variables $\mathbf{x}_2$ and $\mathbf{x}_5$ are doubled to create extra ill-conditioning in design matrix.}
\centering
\fbox{%
\begin{tabular}{*{6}{c}}
\hline
Effect &       Mean         & Variance  &  Effect  &    Mean     & Variance \\ \hline
$\tau_1$   & 0.061179101 & 47.85924599 &$\beta_1$ &    0.2566886 & 7.801541e+02 \\
$\tau_2$   & 2.394966406 & 48.45928806 &$\beta_2$ &    -0.1343304& 1.987533e+02 \\
$\tau_3$   & 0.991546801 & 0.07191629 &$\beta_3$ &    3.4266491  & 3.609550e+03 \\
$\tau_4$   & 1.991347812 & 0.05987755 &$\beta_4$ &    1.9952690  & 2.088842e+03 \\
$\tau_1^w$ & -0.003064839& 0.01757001 &$\beta_5$ &    1.7629811  & 4.647251e+02 \\
$\tau_2^w$ & 2.240328457 & 0.02063512 &$\beta_6$ &    0.9949393  & 2.804560e-01 \\
$\tau_3^w$ & 0.991225900 & 0.06309851 &$\beta_9$ &    1.9752619  & 2.023134e-01 \\
$\tau_4^w$ & 1.978481900 & 0.05714094 &$\beta_{10}$ & 2.9970071  & 9.783558e-02 \\ \hline
\end{tabular}}
\end{table}

Case 5 (Table 6) illustrates the dependence of the weighted group effect on the (within group) APC arrangement of the variables. For this case, we changed the signs of $\mathbf{x}_2$ and $\mathbf{x}_5$, so everything else about the model is still the same except variables in the two strongly correlated groups are not all positively correlated.  The cancellation phenomenon from the uniform model [see (\ref{tp1})] responsible for the optimality of the average and weighted group effects disappears. Both the average and the weighted group effects have large variances, and become not estimable.

\begin{table}
\caption{\label{table6} Case 5: Mean and variance of 8 estimated group effects and 8 estimated individual effects based on 1000 simulated values. Model setting: $(w_1,w_2)=(0.90,0.90)$ and the signs of $\mathbf{x}_2$ and $\mathbf{x}_5$ are switched to create negative correlations in the two correlated groups.}
\centering
\fbox{%
\begin{tabular}{*{6}{c}}
\hline
Effect &       Mean         & Variance  &  Effect  &    Mean     & Variance \\ \hline
$\tau_1$   & 0.2291634 &41.52000724&$\beta_1$ &    0.2148349 &33.9349608 \\
$\tau_2$   & 2.1522303 &9.49975861 &$\beta_2$ &    0.2434918 &50.4173720 \\
$\tau_3$   & 1.0040161 &0.13543209 &$\beta_3$ &    1.1047177 &27.9960202 \\
$\tau_4$   & 2.0056387 &0.60711293 &$\beta_4$ &    2.1359976 &29.4086010 \\
$\tau_1^w$ & 0.2282927 &41.02163441&$\beta_5$ &    3.2159757 &21.6549465 \\
$\tau_2^w$ & 2.1165225 &9.45371635 &$\beta_6$ &    0.9924757 & 0.7200925 \\
$\tau_3^w$ & 1.0047687 &0.09480719 &$\beta_9$ &    2.0068177 & 0.6354657 \\
$\tau_4^w$ & 1.9935899 &0.59861005 &$\beta_{10}$ & 3.0010600 & 0.9066587 \\ \hline
\end{tabular}}
\end{table}

The above examples are a part of a larger numerical study on the variability weighted group effect in a general model (\ref{m0}). The study shows that under the within group APC arrangement, while the estimator of the variability weighted group effect can no longer achieve the limiting variance in Theorem \ref{thm3} for model (\ref{m2}), it is nevertheless remarkably accurate. Its variance is also a decreasing function of the level of multicollinearity. These suggest a perfect uniform correlation structure is only needed for archiving the limiting variance. Without the uniform structure, the weighted group effect is still accurately estimated. The local nature of this type of multicollinearity problem is also clearly shown by the examples; see the Appendix for further discussions.

\subsection{Constrained local regression method}

Ridge regression of Horel and
Kennard (1970) may be used to estimate parameters of strongly correlated variables. However, this method shrinks all estimated values of parameters in the model, whether or not they are that of the strongly correlated variables. This is undesirable as the impact of the multicollinearity problem here is local; least-squares estimates for variables not in the strongly correlated group are unaffected and need not be changed. In the following, we propose a follow-up procedure to the least-squares method which retains the unaffected estimates and deals with the inference and estimation of the affected parameters separately using the weighted group effect.

Suppose there are $p$ predictor variables in a strongly correlated group in its APC arrangement, and denote by $\boldsymbol{\beta}_p^t\in \mathbb{R}^p$ the vector of parameters for variables in this particular group. Then, the variability weighted group effect for this group is
\beq
\tau_w=\mathbf{w}_w^T\boldsymbol{\beta}_p^t,  \label{optline}
\eeq
where $\mathbf{w}_w \in \mathbb{R}^p$ is the known weight vector computed using this group of variables and (\ref{weights}), and $\boldsymbol{\beta}_p^t$ and $\tau_w$ are unknown. Since $\tau_w$ can be accurately estimated by $\hat{\tau}_w$, replacing the unknown $\tau_w$ with the estimated $\hat{\tau}_w$ and the unknown $\boldsymbol{\beta}_p^t$ with a general $\boldsymbol{\beta}_p \in \mathbb{R}^p$ in equation (\ref{optline}), we obtain a line
\beq
\hat{\tau}_w=\mathbf{w}_w^T\boldsymbol{\beta}_p,  \label{optline1}
\eeq
which the true value $\boldsymbol{\beta}_p^t$ should be close to, that is, $\hat{\tau}_w\approx\mathbf{w}_w^T\boldsymbol{\beta}_p^t$. To make use of this information to estimate $\boldsymbol{\beta}_p^t$, we first look for the point $\boldsymbol{\beta}^*_p$ on line (\ref{optline1}) with the shortest distance to the origin. This point serves as a ``lower bound'' for feasible estimates for $\boldsymbol{\beta}_p^t$; we should not shrink the estimate for $\boldsymbol{\beta}_p^t$ so much that its norm is below $\|\boldsymbol{\beta}^*_p\|$. Explicit formula is available for computing $\boldsymbol{\beta}^*_p$ using $\hat{\tau}_w$ and $\mathbf{w}_w$. Once $\boldsymbol{\beta}^*_p$ is computed, we explore other points on line (\ref{optline1}) by setting $\|\boldsymbol{\beta}_p\|^2=c$ where $c$ is a constant larger than $\|\boldsymbol{\beta}^*_p\|^2$. The line in (\ref{optline1}) intersects the $p$-sphere defined by $\|{\boldsymbol{\beta}}_p\|^2=c$ at two points. Each point represents an estimate for the unknown  $\boldsymbol{\beta}^t_p$, and the final selection of $c$ and the point as the estimated value of $\boldsymbol{\beta}_p^t$ can be made through various means such as cross-validation.

Table \ref{table7} is an example that illustrates the above method using the model in the last subsection containing two correlated groups of variables. This table is computed as follows. We first generated a design matrix $\mathbf{X}$ using the model with $(w_1,w_2)=(0.85,0.8)$. The observed correlations among variables in the first two groups of variables range from 0.85 to 0.98, so there are strong correlations among these variables. We then generated one $\mathbf{y}$ observation, computed the least-squares estimates for individual effects as well as the weighted group effects for the two groups in Table 7.

\begin{table}
\caption{\label{table7} An example of constrained local regression. Model setting: $(w_1,w_2)=(0.85,0.80)$. Estimated values of individual and groups effects with two correlated groups of variables.}
\centering
\fbox{%
\begin{tabular}{*{6}{c}}
\hline
Group & Effect& Estimated  & Std error& $t$-value &   $p$-value \\ \hline
 &  $\beta_0$ &  4.7436    & 0.4350   &10.905 &0.000402 \\
1&  $\beta_1$ & -2.6244    & 3.0230   &-0.868 &0.434299 \\
1&  $\beta_2$ &  2.9447    & 3.3085   &0.890  &0.423737 \\
2&  $\beta_3$ &  0.1987    & 2.9019   &0.068  &0.948688 \\
2&  $\beta_4$ &  2.0414    & 3.5087   &0.582  &0.591916 \\
2&  $\beta_5$ &  3.6507    & 3.6986   &0.987  &0.379491 \\
3&  $\beta_6$ &  0.8789    & 0.4321   &2.034  &0.111721 \\
3&  $\beta_7$ &  1.1888    & 0.5195   &2.288  &0.084007 \\
3&  $\beta_8$ &  1.6645    & 0.5065   &3.286  &0.030316 \\
3&  $\beta_9$ &  1.5009    & 0.4081   &3.678  &0.021242 \\
3&  $\beta_{10}$ & 3.0422  & 0.4283   &7.103  &0.002075  \\ \hline
1&  $\tau_1^w$ &  -0.0153   & 0.1676   &0.0915  & 0.931500 \\
2&  $\tau_2^w$ &  1.8511   & 0.1470  &12.590  & 0.000229 \\ \hline
\end{tabular}}
\end{table}

Table \ref{table7} indicates that individuals effects of variables in Groups 1 and 2 are all not significant. For Group 2, this is due to the large variances for these effects caused by the underlying multicollinearity since we know the true values are not zero. At the bottom of the table, we give the weighted group effects for these two groups of correlated variables and their standard errors. These group effects are accurately estimated with small variances. The $t$-test for the first group effect $\tau_1^w$ has a large $p$-value. This suggests two possibilities: either the two underlying parameters of Group 1 are both zero or they are not zero but they cancelled out after weighting. When the second possibility is considered unlikely, we may conclude the parameters of Group 1 are both zero. We do so in this example, so we drop these two variables from the model which makes it unnecessary to estimate their parameters. The second group effect $\tau^w_2$ has a very small $p$-value. This indicates the underlying parameters of Group 2 are not all zero, and we proceed to estimate these parameters.

The weight vector for $\tau_2^w$ is $\mathbf{w}_w=(0.3712, 0.3218, 0.3068)^T$, so
the corresponding estimated weighted group effect line (\ref{optline1}) is
\beq
1.8511 = 0.3712\beta_3+0.3218\beta_4+0.3068\beta_5.  \label{optline2}
\eeq
Plugging the true value of $\boldsymbol{\beta}_p^t=(\beta_3,\beta_4,\beta_5)=(1,2,3)$ into the right-hand side of this equation, we obtain a value of 1.935, so the true value is close to this line by this measure. For line (\ref{optline2}),
\[
\boldsymbol{\beta}^*_p=({\beta}_3^*,{\beta}_4^*,{\beta}_5^*)^T=( 2.047952, 1.775069, 1.692757)^T,
\]
which may be viewed as an initial estimate of $\boldsymbol{\beta}^t_p$. The squared minimum distance of line (\ref{optline2}) to the origin is $\|\boldsymbol{\beta}^*_p\|^2=10.2104$. Exploring other points on the line by setting $c=\|\boldsymbol{\beta}_p^*\|^2+3=13.2104$ (this $c$ is chosen for convenience, it is not optimized), one of the two points we obtain is
\[
\hat{\boldsymbol{\beta}}_p=(\hat{\beta}_3,\hat{\beta}_4,\hat{\beta}_5)^T=( 0.8742301, 1.9232452, 2.9575739)^T,
\]
which would be a good estimate of $\boldsymbol{\beta}_p^t$. We call the above a constrained local regression as it uses the estimated effect line (\ref{optline2}) as a constraint to estimate only the parameters of a strongly correlated group of variables. A detailed study of this method is in progress and results will appear elsewhere.

\section{Concluding remarks}

We note that for a group of $p$ strongly correlated predictor variables in model (\ref{m0}), the exact optimal effect in $\Xi'$, $\tau^*={\xi}(\mathbf{w}^*)$, can be computed numerically. Here, $\Xi'$ is the one defined for parameters of this group only. There are a total of $2^p$ sign arrangements of these $p$ variables. The APC arrangement is an example of these.
We may view $\Xi'$ as the union of $2^p$ sets $\Xi'_k$, $k=1,2,\dots,2^p$, each defined by one sign arrangement of the variables and all with the same constraint $\sum^p_{i=1} w_i=1$ and $w_i\geq 0$. To find $\tau^*$, we compute the optimal effect $\tau^*_j$ in each $\Xi'_k$ and $\tau^*$ is the $\tau^*_j$ whose estimator (\ref{eeff}) has the smallest variance among all $2^p$ such optimal effects.
In numerical examples we have examined, $\tau^*$ is always an optimal from a $\Xi'_k$ with an APC arrangement. When the correlations among the $p$ variables are all strong, $\tau^*$ is approximately the same as the weighted group effect ${\tau}_w={\xi}(\mathbf{w}_w)$ under the same APC arrangement in that $\mathbf{w}^*\approx \mathbf{w}_w$ and $var(\hat{\tau}^*)\approx var(\hat{\tau}_w)$; in this sense we say that $\tau_w$ is approximately optimal. Substantial difference between $\hat{\tau}^*$ and $\hat{\tau}_w$ is observed when the underlying variables are not strongly correlated, but in this case individual parameters are accurately estimated and group effects such as $\hat{\tau}^*$ and $\hat{\tau}_w$ are of little interest.

There does not seem to be a meaningful interpretation for the exact optimal effect $\tau^*$ beyond that it is the numerical optimal. It is in general difficult to compute. Without a simple analytic expression, it also does not lend itself to a theoretical analysis. In contrast, the weighted group effect $\tau_w$ is easy to compute, always accurately estimated and approximately the optimal under the APC arrangement of the underlying variables. Its theoretical properties under the uniform model also add to its appeal.
We recommend it for testing the hypothesis of no group effect for a group of strongly correlated predictor variables and for estimating the parameters of such variables.

Finally, we have taken a modelling approach to study the multicolinearity problem instead of the traditional approach which focuses on the detection and remedies for the problem. This has provided a more detailed look at the type of the problem generated by strongly correlated variables and allowed us to characterize estimable effects centered around the weighted group effect. We are studying this type of multicollinearity problem under other commonly used parametric models for the correlation structure, and hope to further develop this modelling based approach.

\noindent Min Tsao, Department of Mathematics and Statistics, University of Victoria, Victoria, British Columbia, Canada V8W 3R4. Email: mtsao@uvic.ca.

\section{Appendix}

\subsection{Proof of Theorem \ref{thm11}}

For brevity, we only sketch a geometric proof below. For any given $i$, since $\bar{x}_i=0$ and $\|\mathbf{x}_i\|=1$,
\[ corr(\mathbf{x}_i, \mathbf{x}_1)=\mathbf{x}_i \cdot
\mathbf{x}_1=\|\mathbf{x}_i\|\|\mathbf{x}_1\|\cos(\theta_i)=\cos(\theta_i), \]
where $\theta_i$ is the angle between $\mathbf{x}_i$ and $\mathbf{x}_1$. It follows from (\ref{cond1}) that
\beq
|\cos(\theta_i)| > \frac{\sqrt{2}}{2} \label{temp30}
\eeq
for $i=2,3,\dots,p$. Let $\theta$ be the angle between a vector $\mathbf{x}$ and $\mathbf{x}_1$. The region
\[
\left\{\mathbf{x}\hspace{0.05in}|\hspace{0.05in}
\mathbf{x}\in \mathbb{R}^n \mbox{\hspace{0.05in} and \hspace{0.05in}} |\cos(\theta)| > \frac{\sqrt{2}}{2}
\right\}
\]
is the interior of a double cone with apex at the origin, an apex angle of $\pi/4$ and an axis containing
$\mathbf{x}_1$. By (\ref{temp30}), all $\mathbf{x}_i$ are inside this double cone.

Call the cone that contains $\mathbf{x}_1$ the first half of the double cone and its opposing cone the second
half. If an $\mathbf{x}_i$ is in the second half of the double cone, then $-\mathbf{x_i}$ must be in the first
half. Without loss of generality, suppose $\{ \mathbf{x}_2, \mathbf{x}_3, \dots, \mathbf{x}_k\}$ are in the
first half and $\{\mathbf{x}_{k+1}, \mathbf{x}_{k+2}, \dots, \mathbf{x}_p \}$ are in the second half for some
$k$ value satisfying $2\leq k\leq p$. Then
\beq
\{ \mathbf{x}_1, \mathbf{x}_2, \dots, \mathbf{x}_k, -\mathbf{x}_{k+1}, -\mathbf{x}_{k+2}, \dots, -\mathbf{x}_p
\}
\label{temp31}
\eeq
are all in the first half of the double cone. Since the apex angle is only $\pi/4$, the angle $\eta$ between a
pair of variables from the set in (\ref{temp31}) satisfies $0\leq \eta<\pi/2$. This and the fact that they are
all unit vectors imply that their correlation coefficient, which equals $\cos(\eta)$, must be positive. Thus,
the set in (\ref{temp31}) is one of two such sets where all pairwise correlations are positive. \hfill $\Box$

\subsection{Proof of Theorem \ref{thm30}}

Under conditions (\ref{c1}) and (\ref{c2}), matrix $\mathbf{X}^T\mathbf{X}$ in (\ref{xtx}) is
\beq
\mathbf{X}^T\mathbf{X} =
\left( \begin{array}{ccccc}
1 & r & r & \dots & r \\
r & 1 & r & \dots & r  \\
r & r & 1 & \dots & r  \\
\cdot&\cdot&\cdot&\dots &\cdot \\
r & r & r & \dots & 1 \\
\end{array} \right)_{p\times p},  \label{xtx2}
\eeq
where $0<r<1$.
The variance of the least-squares estimator $\hat{\boldsymbol{\beta}}$ in (\ref{ole}) is
\beq
var(\hat{\boldsymbol{\beta}})=(\mathbf{X}^T\mathbf{X})^{-1}. \label{olevar}
\eeq
It follows from (\ref{eeff}) and (\ref{olevar}) that the variance of $\hat{\xi}(\mathbf{w})$ is
\beq
var(\hat{\xi}(\mathbf{w}))=\mathbf{w}^T(\mathbf{X}^T\mathbf{X})^{-1}\mathbf{w}. \label{effvar}
\eeq

In order to see clearly the impact of $r$ on $var(\hat{\xi}(\mathbf{w}))$,
we now express $(\mathbf{X}^T\mathbf{X})^{-1}$ in terms of $r$. By the symmetry imbedded in
$\mathbf{X}^T\mathbf{X}$ in (\ref{xtx2}), its inverse should be of a similar form, that is,
\beq
(\mathbf{X}^T\mathbf{X})^{-1} =
\left( \begin{array}{ccccc}
t & v & v & \dots & v \\
v & t & v & \dots & v  \\
v & v & t & \dots & v  \\
\cdot&\cdot&\cdot&\dots &\cdot \\
v & v & v & \dots & t \\
\end{array} \right)_{p\times p}, \label{invxtx2}
\eeq
where $t$ and $v$ are unknown constants.
Suppose the form in (\ref{invxtx2}) is correct. It follows from
$[\mathbf{X}^T\mathbf{X}][\mathbf{X}^T\mathbf{X}]^{-1}=\mathbf{I}$ that
\bey
& & t+(p-1)rv = 1,  \label{eq1} \\
& & v+tr+(p-2)rv=0.  \label{eq2}
\eey
Solving the system of equations (\ref{eq1}) and (\ref{eq2}) for $t$ and $v$, we obtain
\beq
t=\frac{1+(p-2)r}{1+(p-2)r-(p-1)r^2}, \label{t}
\eeq
and
\beq
v=\frac{-r}{1+(p-2)r-(p-1)r^2}. \label{v}
\eeq
It can be readily verified that the matrix in (\ref{invxtx2}) with elements given by (\ref{t}) and (\ref{v}) is
indeed the unique inverse of $\mathbf{X}^T\mathbf{X}$ in (\ref{xtx2}).
Thus the variance of an estimator given in (\ref{effvar}) can be expressed as
\beq
var(\hat{\xi}(\mathbf{w}))=\mathbf{w}^T(\mathbf{X}^T\mathbf{X})^{-1}\mathbf{w}=t\sum^p_{k=1}w_k^2+v\sum^p_{k=1}\left(w_k\sum_{i\neq k}w_i\right). \nonumber
\eeq
This and (\ref{t}) and (\ref{v}) then imply the theorem. \hfill $\Box$

\subsection{Proofs of Corollaries \ref{coro2} and \ref{coro3}}

\noindent{\em Proof of Corollary \ref{coro2}.} The weight vector associated with $\hat{\beta}_j$ is $\mathbf{w}_j=(w_1,w_2,\dots,w_p)^T$ where $w_j=1$ and $w_k=0$ for $k\neq j$.  For any fixed $r$, at $\mathbf{w}=\mathbf{w}_j$, the $\sum^p_{i=1}w_i^2$ term in the numerator of (\ref{effvar0}) reaches its maximum value of 1 and the $\sum_{i<j}w_iw_j$ term reaches its minimum of 0. Further, the denominator in (\ref{effvar0}) is a positive constant since $0<r<1$. These imply (\ref{re2}).

The above discussion also gives the following simplified version of (\ref{effvar0}) for $var(\hat{\beta}_j,r)$,
\beq
var(\hat{\beta}_j,r) = \frac{1+(p-2)r}{1+(p-2)r-(p-1)r^2}.  \nonumber
\eeq
It can be verified using the above expression that $var(\hat{\beta}_j,r)$ is a strictly monotone increasing function of $r$ with the limit showing in (\ref{limit1}).       \hfill $\Box$

\vspace{0.2in}

\noindent{\em Proof of Corollary \ref{coro3}.}
For any $\mathbf{w}=(w_1, w_2,\dots,w_p)^T$ satisfying $\sum^p_{1=i}w_i=1$ and $w_i\geq 0$, we have
\beq
1^2=\left(\sum^p_{i=1}w_i\right)^2=\sum^p_{i=1}w_i^2+2\sum_{i<j}w_iw_j. \label{t3}
\eeq
Thus, when the first term on the right-hand side $\sum^p_{1=i}w_i^2$ is minimized, the second term $2\sum_{i<j}w_iw_j$ is maximized.
The first term is minimized at $\mathbf{w}=\mathbf{w}_0=(1/p, 1/p,\dots,1/p)^T\in \mathbb{R}^p$. So the numerator of (\ref{effvar0}) is minimized at $\mathbf{w}_0$. This proves (\ref{re1}).

At $\mathbf{w}=\mathbf{w}_0$, the numerator of (\ref{effvar0}) is
\beq
\frac{1+(p-2)r}{p}-\frac{p-1}{p}r=\frac{1-r}{p}. \label{tp1}
\eeq
The denominator of (\ref{effvar0}) can be factored as follows,
\beq
1+(p-2)r-(p-1)r^2=(1-r)[1+(p-1)r].  \label{tp2}
\eeq
It follows from (\ref{effvar0}), (\ref{tp1}) and (\ref{tp2}) that
\beq
var(\hat{\tau}_a,r)=var(\hat{\xi}(\mathbf{w}_0),r)=\frac{1}{p+p(p-1)r}. \label{gev}
\eeq
Equation (\ref{gev}) shows $var(\hat{\tau}_a,r)$ is a monotone decreasing function of $r$ and it also implies (\ref{limit2}).

Finally, let $\mathbf{w}$ be an arbitrary but fixed weight vector such that $\mathbf{w}\neq \mathbf{w}_0$. Since function $\sum^p_{i=1}w_i^2$ attains its minimum of $\sum^p_{i=1}\frac{1}{p^2}$ at only $\mathbf{w}=\mathbf{w}_0$, there exists a small constant $\delta>0$ for this $\mathbf{w}$ such that
\[
\frac{\delta}{p} = \sum^p_{i=1}w_i^2-\sum^p_{i=1}\frac{1}{p^2}=\sum^p_{i=1}w_i^2-\frac{1}{p}>0.
\]
This and (\ref{t3}) imply that for this fixed $\mathbf{w}$,
\[
\sum^p_{i=1}w_i^2=\frac{1+\delta}{p} \mbox{\hspace{0.2in} and \hspace{0.2in}}
2\sum_{i<j}w_iw_j = 1-\frac{1+\delta}{p}.
\]
Thus by (\ref{effvar0}), the variance for the estimator of the effect defined by this $\mathbf{w}$ is
\bey
var(\hat{\xi}(\mathbf{w}),r)&=& \frac{[1+(p-2)r](1+\delta)-r(p-1-\delta)}{p[1+(p-2)r-(p-1)r^2]} \nonumber \\
&=& \frac{(1-r)+\delta+(p-1)r\delta}{p(1-r)[1+(p-1)r]}. \label{t5}
\eey
The numerator of (\ref{t5}) approaches $p\delta>0$ as $r$ approaches 1, but the denominator approaches 0. Thus, for this fixed $\mathbf{w}$ the variance $var(\hat{\xi}(\mathbf{w}),r)$ goes to infinity as $r$ approaches 1. The last part of the corollary follows from this observation, (\ref{limit2}) and the arbitrary selection of $\mathbf{w}\neq \mathbf{w}_0$.  \hfill $\Box$

To see that all estimable effects ${\xi}(\mathbf{w})$ are all around ${\xi}(\mathbf{w}_0)$ in the sense of their associated weight vectors $\mathbf{w}$ are around $\mathbf{w}_0$, by (\ref{t5}) the variance $var(\hat{\xi}(\mathbf{w}),r)$ is a strictly monotone increasing function of $\delta$ which measures the distance between $\mathbf{w}$ and $\mathbf{w}_0$. The collection of estimable effects are given by weight vectors $\mathbf{w}$ with $\delta$ values that are small relative to $(1-r)$ and are thus all around ${\xi}(\mathbf{w}_0)$.

\subsection{Proof of Theorem \ref{thm3}}

Let $\bar{y}$ be the mean of elements of $\mathbf{y}$, and let 
\[\mathbf{y}'=\mathbf{y}-\bar{y}\mathbf{1}_n, \mbox{\hspace{0.1in} and \hspace{0.1in}}
\mathbf{x}'_j=\frac{\mathbf{x}_j-\bar{x}_j\mathbf{1}_n}{s_j}. \]
Then, $\mathbf{y}'$ is the centered response variable with mean 0, and $\mathbf{x}'_j$ the standardized predictor variable with mean 0 and length 1. In terms of these variables, model (\ref{m2}) becomes the
standardized model
\begin{equation}
\mathbf{y}'=\mathbf{X}'\boldsymbol{\beta}' +\boldsymbol{\varepsilon}' \label{m3}
\end{equation}
where $\mathbf{X}'=[\mathbf{x}'_1,\mathbf{x}'_2,\dots,\mathbf{x}'_p]$ is the new $n\times p$ design matrix whose
columns are unit vectors with mean zero. The intercept term of (\ref{m2}) is removed when $\mathbf{y}$ is centered but it can be easily recovered using parameters of  (\ref{m3}). Since $corr(\mathbf{x}'_i,\mathbf{x}'_j)=corr(\mathbf{x}_i,\mathbf{x}_j)=r$ for all $(i,j)$,
$\mathbf{X'}^T\mathbf{X'}$ has the uniform structure in (\ref{xtx2}). We now show that results obtained for model (\ref{m1}) in Section 2.2 apply to the standardized model (\ref{m3}).

Let $\hat{\boldsymbol{\beta}}'$ be the least squares estimator for $\boldsymbol{\beta}'$, i.e.,
\beq
\hat{\boldsymbol{\beta}}'=[\mathbf{X}'^T\mathbf{X}]^{-1}\mathbf{X}'^T\mathbf{y}', \label{extra1}
\eeq
Since $\mathbf{X'}^T\mathbf{X'}$ has the uniform structure, from the proof of Theorem \ref{thm30} we see that results for model (\ref{m1}) in Section 2.2 would apply to the standardized model (\ref{m3}) if
\beq
Var(\hat{\boldsymbol{\beta}}')=\sigma^2[\mathbf{X}'^T\mathbf{X}]^{-1}. \label{extra2}
\eeq
To show this is true, first note that
$$\mathbf{y}'=\mathbf{y}-\bar{y}\mathbf{1}_n=[\mathbf{I}-\frac{1}{n}\mathbbm{1}_{n\times n}]\mathbf{y},$$
where $\mathbbm{1}$ is an $n\times n$ matrix of 1's. So by (\ref{extra1}),
$$ Var(\hat{\boldsymbol{\beta}}')=[\mathbf{X}'^T\mathbf{X}]^{-1}\mathbf{X}'^T [\mathbf{I}-\frac{1}{n}\mathbbm{1}_{n\times n}]Var(y)
[\mathbf{I}-\frac{1}{n}\mathbbm{1}_{n\times n}] \mathbf{X}[\mathbf{X}'^T\mathbf{X}]^{-1}.   $$
Next, noting that $Var(y)=\sigma^2\mathbf{I}$ and that $\mathbbm{1}_{n\times n} \mathbf{X}'=0$ since $\mathbf{X}'$ is the matrix of standardized variables with mean zero, we obtain (\ref{extra2}). Since $\mathbf{X'}^T\mathbf{X'}$ has the uniform structure, (\ref{extra2}) shows that if we add a $\sigma^2$ factor to all expressions of variance in theorems and corollaries for model (\ref{m1}) in Section 2.2, then they apply to model (\ref{m3}).

Parameters $\boldsymbol{\beta}$ and $\boldsymbol{\beta}'$ in models (\ref{m2}) and (\ref{m3}) satisfy $\boldsymbol{\beta}'=\mathbf{S}\boldsymbol{\beta}$
where $\mathbf{S}=diag(s_1,s_2,\dots,s_p)$.
Their least-squares estimators satisfy
\beq
\hat{\boldsymbol{\beta}}'=\mathbf{S}\hat{\boldsymbol{\beta}}. \label{exc2}
\eeq

The minimum-variance unbiased linear estimator for the average group effect of the standardized model (\ref{m3}) is
\beq
\hat{\tau}'_a=\mathbf{w}_0\hat{\boldsymbol{\beta}}'. \label{temp3}
\eeq
Using Corollary \ref{thm30}, we can show its variance $var(\hat{\tau}'_a,r)$ is a strictly monotone decreasing function of
$r$ with limit
\beq
\lim_{r\rightarrow 1} var(\hat{\tau}'_a,r) = \frac{\sigma^2}{p^2}.  \label{temp11}
\eeq

We now express the variance of an estimator $\hat{\xi}(\mathbf{w})$ in (\ref{eeff}) for an effect $\xi(\mathbf{w})$ of the original model (\ref{m2}) in terms of that for the corresponding effect of the standardized model (\ref{m3}).
By (\ref{exc2}), $\hat{\boldsymbol{\beta}}=\mathbf{S}^{-1}\hat{\boldsymbol{\beta}}'$. So
\beq
\hat{\xi}(\mathbf{w})=\mathbf{w}^T\hat{\boldsymbol{\beta}}=\mathbf{w}^T \mathbf{S}^{-1}\hat{\boldsymbol{\beta}}'=h(\mathbf{w}^T \mathbf{S}^{-1})
\frac{\mathbf{w}^T \mathbf{S}^{-1}}{h(\mathbf{w}^T \mathbf{S}^{-1})}\hat{\boldsymbol{\beta}}',  \label{link}
\eeq
where $h(\mathbf{w}^T \mathbf{S}^{-1})$ denotes the sum of the elements of $\mathbf{w}^T \mathbf{S}^{-1}$.
Since the inverse $\mathbf{S}^{-1}=diag(s_1^{-1},s_2^{-1},\dots,s_p^{-1})$, by (\ref{weights}),
\[h(\mathbf{w}_w^T\mathbf{S}^{-1})=h\left(\frac{1}{\sum^p_{i=1}s_i}\mathbf{1}_p^T\right)=\frac{p}{\sum^p_{i=1}s_i},\]
where $\mathbf{1}_p$ is the $p\times 1$ vector of 1's. This and (\ref{link}) imply that
\beq
\hat{\xi}(\mathbf{w}_w)=h(\mathbf{w}^T_w \mathbf{S}^{-1})
\frac{\mathbf{w}^T_w \mathbf{S}^{-1}}{h(\mathbf{w}^T_w \mathbf{S}^{-1})}\hat{\boldsymbol{\beta}}'
=\frac{p}{\sum^p_{i=1}s_i}\mathbf{w}_0^T \hat{\boldsymbol{\beta}}'= \frac{p}{\sum^p_{i=1}s_i} \hat{\tau}'_a.  \label{link2}
\eeq
It follows from this  that
\beq
\hat{\tau}_w=
\frac{p}{\sum^p_{i=1}s_i}\hat{\tau}'_a. \nonumber
\eeq
Thus
\beq
var(\hat{\tau}_w,r)=\left( \frac{p}{\sum^p_{i=1}s_i} \right)^2 var(\hat{\tau}'_a,r). \nonumber
\eeq
Corollary \ref{coro3} and (\ref{temp11}) then imply $var(\hat{\tau}_w,r)$ is a strictly monotone decreasing function of
$r$ and has the limit showing in (\ref{limit5}).

Finally, to see that $\hat{\tau}_w$ is the only one among estimators for the class $\Xi'=\{\xi(\mathbf{w'})\}$ that has a bounded variance when $r$ goes to 1, it suffices to consider only estimators for effects in the subclass $\Xi=\{\xi(\mathbf{w})\}$.
By (\ref{link}), any estimator $\hat{\xi}(\mathbf{w})$ in this subclass can be expressed as
\beq
\hat{\xi}(\mathbf{w})=h(\mathbf{w}^T \mathbf{S}^{-1})\mathbf{w}_+^T\hat{\boldsymbol{\beta}}'
\label{t33}
\eeq
where $\mathbf{w}_+=\mathbf{S}^{-1}\mathbf{w}/h( \mathbf{w}^T\mathbf{S}^{-1})$
which is a proper weight vector. Thus
\beq
var(\hat{\xi}(\mathbf{w}))=h^2(\mathbf{w}^T \mathbf{S}^{-1})var(\mathbf{w}_+^T\hat{\boldsymbol{\beta}}'), \label{t2.1}
\eeq
which implies $var(\hat{\xi}(\mathbf{w}))$ is bounded when $r$ approaches 1 if and only if the corresponding $var(\mathbf{w}_+^T\hat{\boldsymbol{\beta}}')$ is bounded. By Corollary \ref{coro3}, $var(\mathbf{w}_+^T\hat{\boldsymbol{\beta}}')$ is bounded if and only if $\mathbf{w}_+=\mathbf{w}_0$. Since $\mathbf{w}_+=\mathbf{w}_0$ if and only if $\mathbf{w}=\mathbf{w}_w$, it follows that $var(\hat{\xi}(\mathbf{w}_w))$ is the only variance among $var(\hat{\xi}(\mathbf{w}))$ for estimators of effects in $\Xi'$ that remains bounded when $r$ approaches 1. \hfill $\Box$

For a fixed $r$, the variance of an estimator $var(\hat{\xi}(\mathbf{w}),r)$ in (\ref{t2.1}) also depends on $s_i$, so $var(\hat{\xi}(\mathbf{w}_w),r)$ may not be the smallest and consequently $\tau_w$ may not be the exact optimal. But all estimable effects for model (\ref{m2}) are in a neighborhood of $\tau_w$. To see this, we have noted after the proof of Corollary \ref{coro3} that for the simple uniform model (\ref{m1}), estimable effects are in a neighborhood of the average group effect $\tau_a$. Thus, estimable effects for model (\ref{m3}) are in a neighborhood of $\tau'_a$. By (\ref{t33}), an effect for model (\ref{m2}) is a linearly transformed effect for model (\ref{m3}). Transforming the estimable effects of model (\ref{m3}) using (\ref{t33}), we obtain estimable effects for model (\ref{m2}) which form a neighborhood around $\tau_w$.

At a fixed $r$, the most difficult effect to estimate is the individual effect $\beta_j$ of variable $\mathbf{x}_j$ with the smallest variablility $s_j$. To see this, by (\ref{t2.1}), the minimum-variance unbiased linear estimator for a $\beta_i$ has variance
\beq
var(\hat{\xi}(\mathbf{w}_i))=h^2(\mathbf{w}_i^T \mathbf{S}^{-1})var(\mathbf{w}_i^T\hat{\boldsymbol{\beta}}')
=var(\hat{\beta_i}')/s_i^2 .  \label{t2.2}
\eeq
Under the standardized model (\ref{m3}), by Corollary \ref{coro2} and its proof, the $var(\hat{\beta_i}')$ are the same for $i=1,2,\dots,p$ and they are the largest among variances of estimators for effects of model (\ref{m3}). This observation, (\ref{t2.1}) and (\ref{t2.2}) then imply that $var(\hat{\beta}_j)=var(\hat{\xi}(\mathbf{w}_j))$ is the largest among all $var(\hat{\xi}(\mathbf{w}))$ for estimators (\ref{eeff}) of effects in $\Xi'$.

\subsection{The local nature of the multicollinearity problem}

We have seen through numerical examples that the multicollinearity problem generated by a group of strongly correlated variables is a local problem in that its impact is limited to least-squares estimates of parameters for variables in the group. This was a motivating factor behind our approach of isolating and modelling such a multicollinearity problem with the uniform models. If the problem is not local, then local models such as the uniform models would be neither justified nor helpful for understanding the impact of the problem. Here, we provide further discussion on the local nature of the problem.

Consider a simple case of model (\ref{m0}) where $\sigma^2=1$, $p=3$ and
\beq
\mathbf{X}^T\mathbf{X} =
\left( \begin{array}{ccc}
1.0 & 0.9 & 0.1 \\
0.9 & 1.0 & -0.2 \\
0.1 & -0.2 & 1.0  \\
\end{array} \right),  \label{expl1.1}
\eeq
so only $\mathbf{x}_1$ and $\mathbf{x}_2$ are strongly correlated and
variable $\mathbf{x}_3$ is outside of the strongly correlated group. The inverse of $\mathbf{X}^T\mathbf{X}$ is
\beq
(\mathbf{X}^T\mathbf{X})^{-1} =
\left( \begin{array}{ccc}
9.230& -8.846& -2.692 \\
-8.846&  9.519& 2.788 \\
-2.692&  2.788&  1.826 \\
\end{array} \right).  \label{expl1.2}
\eeq
Variances of the estimated individual effects for the strongly correlated $\mathbf{x}_1$ and $\mathbf{x}_2$ are
over 9.0, whereas that for $\mathbf{x}_3$ is less than 2.0. Unusually large values in
$(\mathbf{X}^T\mathbf{X})^{-1}$ are responsible for this substantial difference, and they appear only in
elements $(i,j)$ of $(\mathbf{X}^T\mathbf{X})^{-1}$ where both $i$ and $j$ belong to the two strongly correlated
variables. This always occurs when multicollinearity is caused by strongly correlated variables.
To see this, consider a modified (\ref{m1}) with $\bar{x}_i=0$ and $\|\mathbf{x}_i\|=1$ but only $k$ ($2\leq k<p$) of its variables are strongly correlated. The $k\times k$ block of $\mathbf{X}^T\mathbf{X}$ representing the correlation
matrix of these $k$ variables contains elements whose values are all close to 1. They dominate other off-diagonal
elements of $\mathbf{X}^T\mathbf{X}$ which are small and random correlation coefficients involving other
variables. As such, they make the rows/columns they reside in nearly linearly dependent; see (\ref{expl1.1}) for
an example. Since $\mathbf{X}^T\mathbf{X}$ has one or two more such nearly linearly dependent column/rows than
its submatrices for the minors of elements in the $k\times k$ block, it is ``more singular'' than these
submatrices. Because of this, the determinant of $\mathbf{X}^T\mathbf{X}$ is small in absolute value relative to
the cofactors associated with these minors. By Cramer's rule, large values must appear in the corresponding
$k\times k$ block in $(\mathbf{X}^T\mathbf{X})^{-1}$, making individual effects of strongly correlated variables
difficult to estimate as their least-squares estimators' variances are the diagonal elements of this block. But
this problem does not affect other individual effects ({\em e.g.}, that of $\mathbf{x}_3$) as their associated
submatrices contain the same number of nearly linearly dependent row/columns as $\mathbf{X}^T\mathbf{X}$. This
makes their associated cofactors equally small as the determinant of $\mathbf{X}^T\mathbf{X}$, so their
least-squares estimators' variances are not large.

The above discussion shows that this type of multicollinearity problem is a local problem
in that it affects only the strongly correlated variables even when there are other variables in the model.

\end{document}